


\input amstex
\documentstyle{amsppt}
\magnification=\magstep1
\baselineskip=18pt
\hsize=6truein
\vsize=8truein

\topmatter
\title On the topology of conformally compact Einstein 4-manifolds
\endtitle

\author Sun-Yung A. Chang$^{\ssize 1}$ , Jie Qing$^{\ssize 2}$
 and Paul Yang$^{\ssize 3}$
\endauthor

\leftheadtext{Conformally compact Einstein}
\rightheadtext{Chang, Qing and Yang}

\adjustfootnotemark{1}\footnotetext{Princeton University,
Department of Mathematics, Princeton, NJ 08544-1000, supported by
NSF Grant DMS--0070542.}
\adjustfootnotemark{1}\footnotetext{University of California,
Department of Mathematics, Santa Cruz, CA 95064, supported in part
by NSF Grant DMS--0103160.}
\adjustfootnotemark{1}\footnotetext{Princeton University,
Department of Mathematics, Princeton, NJ  08544-1000, supported by
NSF Grant DMS--0070526.} \abstract In this paper we study the
topology of conformally compact Einstein 4-manifolds. When the
conformal infinity has positive Yamabe invariant and the
renormalized volume is also positive we show that the conformally
compact Einstein 4-manifold will have at most finite fundamental
group. Under the further assumption that the renormalized volume
is relatively large, we conclude that the conformally compact
Einstein 4-manifold is diffeomorphic to $B^4$ and its conformal
infinity is diffeomorphic to $S^3$.
\endabstract

\endtopmatter
\document

\vskip 0.1in
\head 0. Introduction \endhead

Conformally compact Einstein manifolds play an important part in
the AdS/CFT correspondence, a promising new area under extensive
development in string theory [Mc] [GKP] [W]. Mathematically the
study of conformally compact Einstein structures began with the
work of Fefferman and Graham ([FG1]) in which they gave a
procedure to find local conformal invariants and the associated
conformally covariant operators. It has been a very challenging
problem to give general existence theory. Recently M. Anderson
([A2]) has shown the existence of conformally compact Einstein
metrics whose boundary is the 3-sphere with arbitrary conformal
structure in the connected component of the round one, thus
extending the perturbation existence result of Graham and Lee
([GL]). On the other hand, in recent works, using several
conformally covariant differential equations and their ties to
Chern-Gauss-Bonnet formula, two of the coauthors with Gursky
studied the conformal geometry of closed 4-manifolds and obtained
remarkable progress in understanding of the topology of
4-manifolds. In our view, conformally compact Einstein 4-manifolds
provide a platform to study conformal geometry of closed
3-manifolds. In this note we will at least use it to translate the
results in four dimension to 3-manifolds.

An oriented manifold $(X, g)$ with boundary $M$ is a conformally
compact Einstein manifold if there is a complete Einstein metric
$g$ in the interior of $X$, and a smooth defining function $s$ for
the boundary $M=\partial X=\{s=0\}$ so that $(X, s^2 g)$ is a
compact Riemannian manifold with boundary. Since defining
functions are not unique the data $(X, g)$ determines a conformal
structure $(M, [\hat g])$ which is called the conformal infinity
of $(X, g)$. It will be advantageous to find a defining function
$s$ whose associated conformal metric enjoys optimal positivity
property. To this end, one of us first observe ([Q]),

\proclaim{Theorem } Suppose that $(X^{n+1}, g)$ is a conformally
compact Einstein manifold, and that $\hat g$ is a Yamabe metric for
$(M, [\hat g])$. Then there exists a conformal
compactification $(X, u^{-2}g)$, which has a totally geodesic boundary
$M$ and whose scalar curvature $R[u^{-2}g] \geq \frac {n+1}{n-1}
R[\hat g]$.
\endproclaim

An important invariant of the conformally compact Einstein
structure is the renormalized volume (see section one for the
precise definition). The volume of a conformally compact Einstein
manifold is infinite. But an appropriate normalization gives rise
to the new invariants as suggested in the works of Maldacena [Mc],
Witten [W], and Gubser, Klebanov and Polyakov [GKP]. The
renormalization was carried out by Henningson and Skenderis in
[HS] and later also by Graham in [Gr]. Anderson in [A1] observed
that the renormalized volume for a conformally compact Einstein
4-manifold appears in the Chern-Gauss-Bonnet formula
$$
8 \pi^2 \chi(X) = \frac 14\int_X  |\Cal W|^2 dv_g + 6 V, \tag 0.1
$$
where $\Cal W$ is the Weyl curvature and the norm of the Weyl
tensor is given by $|\Cal W|^2 = \Cal W_{ijkl}\Cal W^{ijkl}$;
i.e., the usual definition when $W$ is viewed as a section of
$\otimes^4 T^{*}M^4$ and $V$ is the renormalized volume of a
conformally compact Einstein 4-manifold $(X^4, g)$. It follows
that the renormalized volume is an invariant of the underlying
conformal structure of the conformally compact Einstein manifold.
On the other hand, we recall that for a closed 4-manifold $(Y,g)$,
we have the Chern-Gauss-Bonnet formula
$$
8 \pi^2 \chi(Y) = \int_Y ( \frac {1}{4} |\Cal W|^2 + \sigma_2 (A_g)) dv_g,
\tag 0.2
$$
where $A_g=Rc - \frac {R}{6}g$ is the Weyl Schouten tensor and
$\sigma_2 (A_g) $ is the second symmetric function of the
eigenvalues of $A_g$. Thus the renormalized volume is a
combination of a boundary invariant as well as the Ricci part of
the Chern-Gauss-Bonnet integral. According to the above cited
Theorem, the given compactification $X$ has a totally geodesic
boundary. It  is then natural to consider the doubling manifold
$Y$. In this setting, we will show that the renormalized volume
plays an important role in restricting the topology and geometry
of the conformally compact Einstein manifolds. In particular, we
have

\proclaim{Theorem A} Suppose that $(X,g)$ is a conformally compact
Einstein 4-manifold with its conformal infinity of positive Yamabe
constant and suppose the renormalized volume $V$ satisfies
$$
V > \frac {1}{3}\frac{4\pi ^2}{3} \chi (X), \tag 0.3
$$
Then $X$ is homeomorphic to the 4-ball $B^4$ up to a possible
finite cover.
\endproclaim
This result uses the vanishing theorems of Gursky [Gu2] to show
that the doubling manifold $Y$ is a homology sphere, and the main
result of [CGY1] to show that $Y$ admits a conformal metric of
positive Ricci tensor, and hence has at most finite fundamental
group. The homeomorphism classification theory of Donaldson and
Freedman then implies that $Y$ is homeomorphic to the 4-sphere up
to a finite cover. On the other hand, under a stronger hypothesis,
we can conclude that $Y$ is diffeomorphic to the 4-sphere.

\proclaim{Theorem B} Suppose that $(X, g)$ is a conformally compact
Einstein 4-manifold with conformal infinity of positive Yamabe constant.
Then
$$
V > \frac 12 \ \frac {4\pi^2}{3} \chi (X)
\tag 0.4
$$
implies that $X$ is diffeomorphic to $B^4$ and $M$ is diffeomorphic to
$S^3$.
\endproclaim

These considerations demonstrate that the renormalized volume is
an important conformal invariant of the conformally compact
Einstein structure. The recent work of Fefferman and Graham [FG2]
showed that it may be interpreted as a $Q$ curvature integral
associated to a natural boundary operator. It would be of interest
to understand this invariant.
 One may ask whether the positivity
of the renormalized volume is an intrinsic property of the
conformal structure of the boundary $M$ alone. This is related to
the uniqueness question of the conformally compact Einstein
structure. There is some partial answer to this question in the
article of Anderson [A3] (cf. [Wa]).

The paper is organized as follows: In section 1 we introduce conformally
compact Einstein manifolds and relevant properties. In section 2 we
construct
positive eigenfunctions and present a conformal compactification of a
conformally
compact Einstein manifolds. In section 3 we recall the results in [CGY1]
[CGY2] [Gu1] [Gu2] and prove our main theorems.

\vskip 0.1in
\head 1. Conformally compact Einstein manifolds
\endhead

Let us start with the definition of a conformally compact Einstein
manifold. Suppose $X^{n+1}$ is a compact oriented (n+1)-manifold
with boundary $\partial X = M^n$. A Riemannian metric $g$ in the
interior of $X$ is said to be a $C^{m, \alpha}$ conformally
compact if $\bar g = s^2 g$ extends as a $C^{m,\alpha}$ metric on
$\bar X$, where $s$ is a defining function of the boundary in the
sense that: $s>0$ in $X$, $s=0$ and $ds \neq 0$ on $M$. Clearly
defining functions are not unique. For a given defining function,
the metric $\bar g$ restricted to $TM$ induces a metric $\hat g$
on $M$. $\hat g$ rescales upon changing the defining function $r$,
therefore defines a conformal structure on $M$. We call $(M, [\hat
g])$ the conformal infinity of the conformally compact manifold
$(X, g)$. Therefore $(X^{n+1}, g)$ is a conformally compact
Einstein manifold if it is conformally compact and $\text{Ric}(g)
= -ng$. A good reference for basic properties of conformally
compact Einstein manifolds is the paper of Graham [Gr].

The boundary regularity of the conformally compact Einstein metric
is an important issue. Thanks to M. Anderson [A2], in 4 dimension,
we know that $C^{2, \alpha}$ would imply the full smoothness of
the conformally compact Einstein metric. In other words, the
conformally compact Einstein 4-manifolds we will be focus on in
the last section are, always smooth when the conformal infinity is
smooth, therefore, the local expansions of Einstein metrics near
the infinity are always available. Sufficient boundary regularity
is assumed in higher dimension to ensure the expansion though.

Solving a nonlinear first order PDE by the method of characteristics near
the
boundary, one has (cf. Lemma 2.1 and Lemma 2.2 in [Gr]):

\proclaim{Lemma 1.1} Given any conformally compact manifold
$(X, g)$, a metric $\hat g$ in the conformal class $[\hat g]$
determines a unique defining function $s$ in a neighborhood of $M$
such that
$$
g = s^{-2} ((ds)^2 + g_s),
\tag 1.1
$$
where, when $n$ is odd,
$$
g_s = \hat g + g^{(2)} s^2 + \text{even powers of $s$}
+ g^{(n-1)}s^{n-1} +
g^{(n)}s^n + \cdots;
\tag 1.2
$$
$g^{(i)}$ for $i < n$ are determined by local geometry of $(M, \hat g)$
and
$g^{(n)}$ is trace-free, when $n$ is even,
$$
g_s = \hat g + g^{(2)} s^2 + \text{even powers of $s$} + h s^n \log s
+ g^{(n)}s^n + \cdots,
\tag 1.3
$$
$g^{(i)}$ for $i < n$, $h$ and the trace of $g^{(n)}$ are determined
by local geometry of $(M, \hat g)$, and $h$ is trace-free.
\endproclaim

For example, one can calculate that
$$
g^{(2)} = - \frac 1{n-2}(\hat R_{ij} - \frac {\hat R}{2(n-1)}\hat
g_{ij}).
\tag 1.4
$$

As a realization of the so-called AdS/CFT correspondence,
Henningson and Skenderis in [HS] and Graham in [Gr] defined and
calculated the renormalized volume for a conformally compact
Einstein manifold. Namely, they considered the expansions, when
$n$ is odd,
$$
\text{Vol}(\{ s > \epsilon\})
= c_0 \epsilon^{-n} + c_2 \epsilon^{-n+2} +
\ \text{odd powers of $\epsilon$} \ + c_{n-1} \epsilon^{-1} + V + o(1),
\tag 1.5
$$
when $n$ is even,
$$
\text{Vol}(\{ s > \epsilon\})
= c_0 \epsilon^{-n} + c_2 \epsilon^{-n+2} +
\ \text{even powers of $\epsilon$} \
+ c_{n-2}\epsilon^{-2} + L\log \frac 1\epsilon + V + o(1),
\tag 1.6
$$
where
$$
c_i = \int_M v^{(i)} dv_{\hat g},
\tag 1.7
$$
and $v^{(i)}$ are local scalar invariants of $(M, \hat g)$ of order $i$.
More importantly he showed the renormalized volume $V$ for $n$
odd and $L$ for $n$ even are independent of the choice of defining
functions.
Incidentally when $n$ is even the quantity $L = \int_M v^{(n)}dv_{\hat
g}$
is a conformal invariant of $M$, and when $n$ is odd, $V$ is only an
invariant
of the ambient structure $(X,g)$.

Particularly, when $n=3$, Anderson in [A]
relates the renormalized volume to the Euler number of
a conformally  compact Einstein 4-manifold,
$$
8 \pi^2 \chi(X) = \frac 14 \int_X  |\Cal W|^2 dv_g + 6 V. \tag 1.8
$$
\vskip 0.1in
\head 2. Conformal compactifications
\endhead

In this section we present a conformal compactification of a general
conformally compact Einstein manifold.
As observed in [Q], when one considers the hemisphere
compactification of the hyperbolic space, the conformal factor is
a positive eigenfunction. Namely, recall
$$
(H^{n+1}, g_H) = (B^{n+1}, (\frac 2{1-|y|^2})^2|dy|^2),
$$
and
$$
(S^{n+1}_+, g_S) = (B^{n+1}, (\frac 2{1+|y|^2})^2|dy|^2),
$$
that is $g_S = t^{-2}g_H$, where we denote by
$$
t = \frac {1+|y|^2}{1-|y|^2}.
$$
One finds after some calculation (cf. [Q]),
$$
\Delta_{g_H} t = (n+1) t.
\tag 2.1
$$
Given a general conformally compact Einstein manifold, using
the theory of uniformly degenerate elliptic linear PDE developed in
[M] [GL] [L], one can construct positive eigenfunctions. A good reference
for our discussions here is a paper of Lee [L]. The
following lemma is a slight improvement of Lemma 5.2 in [L].

\proclaim{Lemma 2.1} Suppose that $(X^{n+1},g)$ is a conformally
compact Einstein
manifold. Then there exists a positive function $u$ such that
$$
\Delta u = (n+1) u \ \ \text{in} \ X.
\tag 2.2
$$
In addition, when $n$ is odd,
$$
u = \frac 1s + \frac {\hat R}{4n(n-1)}s + w^{(4)} s^3 + \
\text{odd powers of $s$} + w^{(n-1)}s^{n-2} + O (s^n) \tag 2.3
$$
and when $n$ is even,
$$
u = \frac 1s + \frac {\hat R}{4n(n-1)}s + w^{(4)} s^3 + \
\text{odd powers of $s$} + w^{(n)}s^{n-1} + O(s^n), \tag 2.4
$$
where $w^{(i)}$ are all local invariants of Riemannian geometry
of $(M, \hat g)$ of order $i$.
\endproclaim

\demo{Proof} To apply the theory of uniformly degenerate elliptic PDE
(cf. [L], for example) we need to calculate
$$
\Delta \frac 1s
= \frac {s^{n+1}}{\sqrt{\det g_s}}\partial_s (
s^{1-n} \sqrt{\det g_s}\partial_s \frac 1s)
= \frac {n+1}s - \frac 12 \text{Tr}_{\hat g}\partial_s g_s.
$$
Then from the expansion of $g_s$ obtained
by Graham as stated in Lemma 1.1 in the previous section,
when $n$ is odd,
$$
\Delta \frac 1s = \frac {n+1}s + p^{(2)}s + \ \text{odd powers of $s$} \
+ p^{(n-1)} s^{n-2} + O(s^n),
$$
and when $n$ is even,
$$
\Delta \frac 1s = \frac {n+1}s + p^{(2)}s + \ \text{odd powers of $s$} \
+ p^{(n)} s^{n-1} + O(s^n),
$$
where $p^{(k)} = -\frac k2 \text{Tr}_{\hat g}g^{(k)}$ are local
invariants
of Riemannian geometry of $(M, \hat g)$ of order $k$, and
$$
p^{(2)} = \frac {\hat R}{2(n-1)},
$$
for example. Now let us calculate the key ingredient in the theory of
uniformly degenerate elliptic PDE: the indicial equation. For that
purpose we first calculate with separating variables for any smooth function $\phi$ defined on $(M, \hat g)$,
$$
\Delta (\phi s^k)
= \frac {\phi s^{n+1}}{\sqrt{\det g_s}} \partial_s
(s^{1-n}\sqrt{\det g_s} \partial_s s^k)
+ \frac {s^{k+2}}{\sqrt{\det g_s}}\partial_\alpha(
\sqrt{\det g_s} g^{\alpha\beta}_s\partial_\beta \phi),
$$
where
$$
\frac {\phi s^{n+1}}{\sqrt{\det g_s}} \partial_s
(s^{1-n}\sqrt{\det g_s} \partial_s s^k)
= k(k-n) \phi s^k + \frac k2 \phi s^{k+1}
\text{Tr}_{\hat g}\partial_s g_s
$$
and
$$
\frac {s^{k+2}}{\sqrt{\det g_s}}  \partial_\alpha(
\sqrt{\det g_s} g^{\alpha\beta}_s\partial_\beta \phi)
= s^{k+2} \tilde\Delta \phi + s^{k+2}(
\cdots \text{even powers of $s$} \
\cdots ),
$$
$\phi$ is independent of $s$, and $\tilde\Delta$ is the Laplacian
of $(M, \hat g)$. Therefore the indicial equation for us is
$$
(\Delta - (n+1)) (\phi s^k) = (k(k-n)- (n+1)) \phi s^k +
s^{k+2}( \cdots \ \text{even powers of $s$} \ \cdots)
$$
Then we may write, when $n$ is odd
$$
(\Delta - (n+1))(\frac 1s + w^{(2)} s + \ \text{odd
powers of $s$} \ + w^{(n-1)} s^{n-2}) = O(s^n)
$$
and when $n$ is even
$$
(\Delta - (n+1))(\frac 1s + w^{(2)}s + \ \text{odd
powers of $s$} \ + w^{(n)} s^{n-1}) = O(s^n).
$$
Notice that $w^{(i)}$ are all local invariants
of Riemannian geometry of $(M, \hat g)$, for example for $ i =2$ we have
$$
w^{(2)} =  \frac {\hat R}{4n(n-1)}.
$$
Then, the lemma follows from Proposition 3.3 in the paper
of Lee [L]. Note that positivity of $u$ is a simple consequence of the
maximum
principle.
\enddemo

We remark that for the standard hyperbolic space we have
$$
u = t = \frac 1s + \frac 14 s.
$$

We can now calculate the scalar curvature for the metric $u^{-2}g$ on $X$
as in [Q]. We find, for $u$ obtained in Lemma 2.1,
$$
- \Delta u^{-\frac {n-1}2} - \frac {n^2-1}4 u^{-\frac {n-1}2} =
\frac {n^2-1}4 (u^2 - |du|^2) u^{-\frac {n+3}2}.
\tag 2.5
$$
Thus the scalar curvature of $u^{-2}g$ is
$$
R [u^{-2}g] = n(n+1)(u^2 - |du|^2).
\tag 2.6
$$
Notice that $u$ is determined by the choice of a defining function,
therefore
the choice of a metric of the conformal infinity by Lemma 1.1. So we may
choose a Yamabe metric $\hat g$ for $(M, [\hat g])$ in the following
theorem.

\proclaim{Theorem 2.2} Suppose that $(X^{n+1},g)$ is a conformally
compact
Einstein
manifold, and that $u$ is the eigenfunction obtained in Lemma 2.1 for
a Yamabe metric $\hat g$ of the conformal infinity $(M, [\hat g])$.
Then $(X^{n+1}, u^{-2}g)$ is a compact manifold with totally geodesic
boundary $M$ and the scalar curvature is greater than or equal to
$\frac {n+1}{n-1}\hat R$.
\endproclaim

\demo{Proof} First by Lemma 2.1, one can estimate
$$
\aligned
u^2 & - |du|^2 = \frac 1{s^2} + \frac {\hat R}{2n(n-1)} + O(s^2)
- s^2 (\partial_s u)^2 - s^2 g^{\alpha\beta}_s u_\alpha u_\beta \\
& = \frac 1{s^2} + \frac {\hat R}{2n(n-1)} -
\frac 1{s^2} + \frac {\hat R}{2n(n-1)}
+ O(s^2) \\
& = \frac {\hat R}{n(n-1)} + O(s^2).
\endaligned
\tag 2.7
$$
Thus in light of (2.6),
the theorem follows from a maximum principle and
the following Bochner formula,
$$
- \Delta (u^2- |du|^2) = 2|Ddu - ug|^2,
\tag 2.8
$$
which has also been used in [L]. The smoothness assertion follows from the
expansion (2.3).
\enddemo

\vskip 0.1in
\head 3. Topology of conformally compact Einstein 4-manifolds
\endhead
In this section we apply the vanishing theorems of Gursky and the
deformation theory of [CGY] of conformal metrics on the doubling
4-manifold $Y$ in order to draw conclusions about the topology of
the conformally compact manifold $X$. Under increasingly stringent
conditions on the renormalized volume $V$, we demonstrate the
vanishing of $H^1(Y)$, $H^2(Y)$. This is a reflection of the
increasing constraint placed on the Weyl-Schouten tensor of the
compactified structure.

According to the compactification obtained in Theorem 2.2 we have a
compact 4-manifold $(X, \tilde g)$ with totally geodesic boundary
$M$. In addition  the  scalar curvatures satisfies $R[\tilde g]
\geq 2 R[\hat g]$, where $\hat g$ is the Yamabe metric for the
conformal infinity $(M, [\hat g])$. Recall the Chern-Gauss-Bonnet
formula for 4-manifolds with totally geodesic boundary
$$
8 \pi^2 \chi(X) = \int_X  ( \frac 14  |\tilde \Cal W|^2 + \sigma_2
(A_{ \tilde g}) ) dv_{\tilde g} \tag 3.1
$$
where
$$
\sigma_2(A) = \frac 1{24} R^2 - \frac 12 |E|^2
$$
and $E$ is the traceless Ricci curvature. It follows from the conformal
invariance
$$
|\tilde \Cal W|^2 dv_{\tilde g} = |\Cal W|^2 dv_g \tag 3.2
$$
and Anderson's formula (1.8) that
$$
\int_X \sigma_2(A_{\tilde g}) dv_{\tilde g} = 6 V.
\tag 3.3
$$
Consider the doubling manifold $(Y, \bar g)$ of $(X, \tilde g)$. The metric
$\bar g$
on $Y$ belongs to $C^{2,1}$, i.e. its second derivative is of Lipschitz class,
therefore the
curvature tensor is of Lipschitz class. Consequently we have

\proclaim{Proposition 3.1} If the conformal infinity has positive
Yamabe invariant, then the renormalized volume satisfies
$$
V \leq \frac 43 \pi^2.
\tag 3.4
$$
Equality in (3.4) holds if and only if $(X, g)$ is the hyperbolic space,
therefore $(M, [\hat g])$ is the round sphere.
\endproclaim
\demo{Proof} According to  Gursky ([Gu1])
if a 4-manifold $(Y, \bar g )$ is of positive Yamabe constant, then
$$
\int_Y \sigma_2 (A_{ \bar g}) dv_{\bar g} \leq 16\pi^2,
\tag 3.5
$$
and equality holds if and only if $(Y, \bar g) $ is the round 4-sphere.
In the situation at hand,
the doubling space $Y$ has a metric of positive scalar curvature
according
to Theorem 2.2, therefore,
$$
V \leq \frac 1{6} \int_X \sigma_2 (\tilde g) dv_{\tilde g} \leq
\frac 1{6} 8\pi^2
= \frac 43 \pi^2
\tag 3.6
$$
and equality holds if and only if
$$
\int_X |\Cal W|^2 dv_{g} =0.
$$
In this case the simply connected hyperbolic manifold is the hyperbolic
space, and
the conformal infinity is the round sphere. The proof is complete.
\enddemo

Let us assume that $(X,g)$ is a conformally compact Einstein
4-manifold whose conformal infinity $(M, [\hat g])$ has positive
Yamabe invariant. According to the spectral result of Lee in [L]
and a vanishing result of Wang [Wa](cf. [CG] [WY]), there are no
harmonic one forms in $L^2$. Therefore it follows from the
isomorphism established by Mazzeo [M], the space of $L^2$ harmonic
one forms is isomorphic to the relative cohomology $H^1(X,M)$,
hence the latter has to vanish.

We have the long exact sequence in cohomology for the
pair $(X,M)$:
$$
\CD
... \  H^k(X,M) @>j^*>> H^k(X) @>i^*>> H^k(M) @>\delta ^*>>
H^{k+1}(X,M) \ ...
\endCD
\tag 3.7
$$
and the Mayer-Vietoris sequence for the double $Y=X \cup X'$
where $X'$ is the second copy of $X$:
$$
\CD ... \ H^{k-1}(M) @>\gamma>> H^k(Y) @>\psi ^*>> H^k(X) \oplus
H^k(X') @>\phi ^*
>> H^k(M) \ ...
\endCD
\tag 3.8
$$
Since $H^1(X,M)=0$, it follows from the long exact sequence
that
$i^*H^1(X)=H^1(X)$.
The Mayer-Vietoris sequence then implies that
$H^1(Y)=H^1(X)$.

We wish to impose a second positivity assumption on the
conformally compact Einstein structure in addition to the
positivity of Yamabe invariant of the boundary: the renormalized
volume $V$ should be positive. We remark that this implies in view
of the formula (0.1) that $\chi (X)>0$. In fact we have
a stronger consequence according to the first vanishing theorem of
Gursky [Gu2]:

\proclaim{Theorem 3.2} [Gu2] If $(N, g)$ is a closed 4-manifold
with positive Yamabe invariant, and
$$
\int_N \sigma_2(A_{\bar g})dv_{\bar g} >0, \tag 3.9
$$
then $H^1(N)=0$.
\endproclaim

In order to apply the vanishing theorem, we need to assure
ourselves that the doubling metric is smooth. Since the metric on
$X$ is of the class $C^{2,1}$ and the boundary is totally
geodesic,  we may (if necessary) smooth out the metric by a small
perturbation without changing the two positivity assumptions:
positive Yamabe invariant and $\int_Y \sigma_2(A_{\bar g})
dv_{\bar g} > 0$. We conclude thus $H^1(X)=H^1(Y)=0$.

Next we formulate a stronger positivity condition in order
to control the second homology.
Recall for a closed 4-manifold $(Y, \bar g)$, we have
the Chern-Gauss-Bonnet formula (0, 2).
On the other hand we also have the Hirzebruch signature formula
$$
12 \pi^2 \tau (Y) =  \int_Y \frac 14 (|\bar \Cal W^+|^2 - |\bar \Cal
W^-|^2)dv_{\bar g} \tag 3.10
$$
Combining the two formulae we have
$$
4\pi^2(2\chi(Y) + 3 \tau(Y)) = \frac 14 \int_Y |\bar \Cal W^+|^2dv_{\bar g} +
\int_Y \sigma_2(A_{\bar g}) dv_{\bar g} \tag 3.11
$$
and
$$
4\pi^2(2\chi(Y) - 3\tau(Y)) = \frac 14 \int_Y  |\bar \Cal W^-|^2dv_{\bar g} +
\int_Y \sigma_2(A_{\bar g}) dv_{\bar g}. \tag 3.12
$$

\proclaim{Theorem 3.3} [Gu2]  Suppose that $(Y, \bar g)$ is a closed oriented
4-manifold with positive Yamabe constant. Then
$$
\frac 14 \int_Y |\bar \Cal W^+|^2dv_{\bar g} < \int_Y \sigma_2 (A_{\bar g}) dv_{\bar g} \tag 3.13
$$
implies that the self-dual part of the second homology vanishes;
similarly,
$$
\frac 14 \int_Y |\bar \Cal W^-|^2dv_{\bar g} < \int_Y \sigma_2 (A_{\bar g}) dv_{\bar g} \tag 3.14
$$
implies that the anti-self-dual part of the second homology vanishes.
\endproclaim

\proclaim{Proposition 3.4}
Assume that $(X,g)$ is a conformally compact Einstein 4 manifold whose
conformal infinity has positive Yamabe invariant and that
$$
V > \frac 13 \frac {4 \pi^2}{3} \chi (X)
\tag 0.3
$$
then  $V >0$ and $H^2(X, \Bbb R)=0$.
\endproclaim
\demo{Proof}
According to  (1.8),
$$
\frac 14 \int_X |\tilde \Cal W|^2 dv_{\tilde g} = 8\pi^2 \chi (X)
- 6 V  < 12 V.
$$
Therefore $V > 0$. If we translate the assumption (0.3) into assumption
on the doubling manifold $(Y, \bar g)$, we have
$$
\int_Y \sigma_2 (A_{\bar g}) dv_{\bar g} > \frac 83 \pi^2 \chi (Y)
\tag 3.15
$$
here we note that $\chi (Y) = 2 \chi (X)$. It follows from
(0.2) that (3.15) is equivalent to
$$
\frac 14 \int_Y |\bar \Cal W|^2 dv_{\bar g} < 2 \int_Y \sigma_2 (A_{\bar g}) dv_{\bar g}.
\tag 3.16
$$
At this point, we remark once more that we may smooth the doubling
metric if necessary while preserve the condition (3.16).

We can now assume without loss of generality that
$$
\frac 14 \int_Y |\bar \Cal W^+|^2 dv_{\bar g} <  \int_Y \sigma_2
(A_{\bar g}) dv_{\bar g}. \tag 3.17
$$
and apply the vanishing theorem Theorem 3.2 of Gursky, and
conclude that $H^2_+(Y)$ $= 0$. From (3.11) we also have
$$
4\pi^2 (2\chi (Y) + 3 \tau(Y)) \geq \int_Y \sigma_2(A_{\bar g}) dv_{\bar g}
> \frac 83 \pi^2 \chi (Y).
\tag 3.18
$$
A short computation with the two long exact sequences (3.7) and (3.8)
and using the fact $H^1(X)=H^1(Y)=0$ will yield
that $H^2(Y)=H^2(X) \oplus H^2(X')$ and hence $dim H^2(Y)$ is even.
Hence we we may assume $dim H^2_-(Y)=2k$. Then (3.18) yields
$$
2(2+2k) + 3(-2k) > \frac 23 (2+2k)
\tag 3.19
$$
which is impossible for any positive integer $k$. Therefore $ k =0$ and
$H^2 (X, \Bbb R) = H^2 (Y, \Bbb R) =0$ and we have finished the proof of the proposition.
\enddemo

So far we have applied the vanishing theorems of Gursky to conformally
compact Einstein manifolds. Next we develop the implication
the
recent work of Chang-Gursky-Yang [CGY1] [CGY2] on
the doubling manifold $Y$.

\proclaim{Theorem 3.5} [CGY1]
Suppose that $(Y, \bar g)$ is an oriented closed 4-manifold with
positive Yamabe constant and that
$$
\int_Y \sigma_2(A_{\bar g}) dv_{\bar g} > 0.
$$
Then there exists a conformal metric $g'=e^{2w} \bar g$ with positive
Ricci tensor.
\endproclaim

\vskip .2in

\noindent{\bf Proof of Theorem A}

As a consequence of Theorem 3.5, we can conclude that the doubling
manifold $Y$ has a finite fundamental group. If we consider its
simply connected covering $\tilde Y$, the homeomorphic
classification theory of Donaldson and Freedman for simply
connected 4-manifolds implies that the manifold $\tilde Y$ is
determined by its intersection form. It follows from Proposition
3.4 that $H^2 (Y, \Bbb R) = 0$ under the assumption (0.3) of
Theorem A. The vanishing of $H^2(Y, \Bbb R)$ in turn implies that
$\tilde Y$ is homeomorphic to $S^4$. We have thus established
Theorem A.

\vskip .1in
We will now discuss the condition (0.4) in Theorem B which is
stronger than the condition (0.3) in the statement of Theorem A.
To see the implication of (0.4), we first recall a main result in
[CGY2].

\proclaim{Theorem 3.6} [CGY2]
Suppose that $(Y, \bar g)$ is an oriented closed 4-manifold with
positive Yamabe constant and that
$$
\frac 14 \int_Y |\bar \Cal W|^2 dv_{\bar g} < \int_Y \sigma_2(A_{\bar g})
dv_{\bar g}.
$$
Then $Y$ is diffeomorphic to $S^4$.
\endproclaim

\vskip .1in

\noindent
{\bf Proof of Theorem B}

Again we will use the doubling manifold $(Y, \bar g)$. We
know
by our assumption (0.4) that
$(Y, \bar g)$ is a $C^{2,1}$ closed doubling 4-manifold with positive
scalar
curvature and
$$
\frac 14 \int_Y |\bar \Cal W|^2 dv_{\bar g} < \int_Y \sigma_2
(A_{\bar g}) dv_{\bar g}. \tag 3.20
$$
We observe that the smoothness condition $\bar g \in C^{2,1}$
means that we can find a nearby smooth metric $g_0$ of $\bar g$ such
that the above mentioned three properties all preserved.
Now we are ready to apply the argument in
[CGY1] [CGY2]
to deform $g_0$ into the metric of the standard 4-sphere with the doubling
property
preserved all the way. For this purpose we want to sketch the ideas given
in
[CGY1] [CGY2] in the following.

First, one considers the following functional $F: W^{2,2}(Y)
\longrightarrow \Bbb R$:
$$
F[\omega] = \tilde \gamma_1 \tilde I[\omega] + \gamma_1 I[\omega]
+ \gamma_2 II[\omega] + \gamma_3III[\omega].
\tag 3.21
$$
For our situation we need
to consider an even nowhere-vanishing symmetric $(0,2)$ tensor $\eta$
which
may be taken as the metric $g_0$ for example. For particular choice
of the
coefficients
$$
\aligned
\tilde \gamma_1 & = - \frac {\int_Y \sigma_2 (A_0) dv_0
- \int_Y |\Cal W_0|^2 dv_0}{ 2\int_Y |\eta_0|^2 dv_0} < 0, \\
\gamma_1 & = - \frac 18, \\
\gamma_2 & = 1, \\
\gamma_3 & = \frac {3\delta -2}{24} > 0,
\endaligned
\tag 3.22
$$
following the work of [CY], it is proved in [CGY1] [CGY2] that there
exists
smooth function $\omega$ which achieves the infimum of $F[\omega]$ over
$W^{2,2}(Y)$ and satisfies
$$
\sigma_2(A) - |\Cal W|^2 = \frac \delta 4\Delta R - \tilde\gamma_1
|\eta|^2 \tag 3.23
$$
for all $\delta > \frac 23$, here we need to emphasize that terms
in (3.23) are taken with respect to the metric $e^{2\omega}g_0$.
Second, by a priori estimate and non-degeneracy of the linearization of
(3.23)
(cf. [CGY1] [CGY2]) one may conclude that there exists metric
$e^{2\omega_\delta}g_0$ which satisfies (3.23) for all $\delta > 0$.
Third, take a metric $e^{2\omega_\delta} g_0$
obtained by the previous steps for sufficient small $\delta$
and deform it by Yamabe flow in the conformal class. It was proved in
[CGY1] [CGY2]
that, after some time, the metrics along the Yamabe flow satisfy
$$
\sigma_2 (A) - \frac 14  |\Cal W|^2 > 0. \tag 3.24
$$
Finally, one applies a weak pinching result of Margerin [Ma]. Namely, with
(3.24),
one can deform again the metric by Ricci flow into the standard metric
for the round
4-sphere (cf. Section 2 in [CGY2]).

Finally we need to make sure that in each of the above steps the doubling
symmetry is preserved.
For this purpose we simply make sure that the solution
in each step is unique, at least locally. For the first step, one
may use the uniqueness of the extremal of $F$, because of $\kappa = 0$
(cf. Theorem 2.1 in [CY]), to conclude that the extremal of $F$ has to be
an even function
on the doubling manifold $(Y, g_0)$. In second step, we have a local
uniqueness as the
linearized equation is non-degenerate, which also assures the solution
$\omega_\delta$
obtained in the argument in [CGY1] [CGY2] is even for all $\delta > 0$.
The last
two steps are clear since the uniqueness of the two geometric
flows guarantees
that all metrics along the flows respect the doubling symmetry
 if the initial
metrics are.
Thus the proof is completed.

\bigskip
\noindent
{\bf References}:
\roster

\vskip 0.1in
\item"{[A1]}" M. Anderson, $L^2$ curvature and volume renormalization of
the AHE
metrics on 4-manifolds, Math. Res. Lett., 8 (2001) 171-188.

\vskip 0.1in
\item"{[A2]}" M. Anderson, Boundary regularity, uniqueness and non-uniqueness for AH Einstein metrics on 4-manifolds, preprint.

\vskip 0.1in
\item"{[A3]}" M. Anderson, Einstein metrics with prescribed conformal infinity on 4-manifolds, preprint.

\vskip 0.1in
\item"{[CG]}" M. Cai and G.J. Galloway, Boundary of zero scalar curvature
in the
\newline
AdS/CFT correspondence, Preprint hep-th/0003046.

\vskip 0.1in
\item"{[CY]}" S.-Y. A. Chang and P. Yang, Extremal metrics of zeta functional
determinants
on 4-manifolds, Ann. of Math. 142(1995) 171-212.

\vskip 0.1in
\item"{[CGY1]}" S.-Y. A. Chang, M. Gursky and P. Yang, An equation of
Monge-Ampere type in
conformal geometry and 4-manifolds of positive Ricci curvature, Annals of Math. 155 (2002), 709-787.

\vskip 0.1in
\item"{[CGY2]}" S.-Y. A. Chang, M. Gursky and P. Yang, A conformally invariant
Sphere theorem in four dimension, Preprint 2002.

\vskip 0.1in
\item"{[FG1]}" C. Fefferman, and C.R. Graham, Conformal invariants, in
{\it The mathematical heritage of Elie Cartan}, Asterisque, 1985, 95-116.

\vskip 0.1in
\item"{[FG1]}" C. Fefferman, and C.R. Graham, Q-curvature and Poincare
metrics, Math. Res. Lett.,  9  (2002),  no. 2-3, 139--151.

\vskip 0.1in
\item"{[Gr]}" C. R. Graham, Volume and Area renormalizations for
conformally compact Einstein metrics. The Proceedings of the 19th Winter
School
"Geometry and Physics" (Srn\'{i}, 1999).
Rend. Circ. Mat. Palermo (2) Suppl. No. 63 (2000), 31--42.

\vskip 0.1in
\item"{[GL]}" C.R. Graham and J. Lee, Einstein metrics with prescribed
conformal
infinity on the ball. Adv. Math. 87 (1991), no. 2, 186--225.

\vskip 0.1in
\item"{[GKP]}" S.S.Gubser, I.R.Klebanov and A.M.Polyakov, Gauge
theory correlators from non-critical string theory, Phys Lett. B
428 (1998) 105-114, hep-th/9802109

\vskip 0.1in
\item"{[Gu1]}" M. Gursky, The principal eigenvalue of a conformally
invariant
differential operator, with application to semilinear elliptic PDE, Comm.
Math. Phys.
207 (1999), 131-143.

\vskip 0.1in
\item"{[Gu2]}" M. Gursky, The Weyl functional, de Rham cohomology, and
K\"{a}hler-Einstein metrics. Ann. of Math. (2) 148 (1998), no. 1,
315--337.

\vskip 0.1in
\item"{[HS]}" M. Henningson and K. Skenderis, The holographic Weyl
anomaly, J. High. Energy Phys. 07 (1998) 023 hep-th/9806087;
Holography and the Weyl anomaly hep-th/9812032

\vskip 0.1in
\item"{[L]}" J. Lee,
The spectrum of an asymptotically hyperbolic Einstein manifold. Comm.
Anal. Geom. 3 (1995), no. 1-2, 253--271.

\vskip 0.1in
\item"{[Mc]}" J. Maldacena, The large N limit of superconformal
field theories and supergravity, Adv. Theo. Math. Phy. 2 (1998)
231-252, hep-th/9711200

\vskip 0.1in
\item"{[Ma]}" C. Margerin, A sharp characterization of the smooth
4-sphere in curvature forms; CAG, 6 (1998), no.1, 21-65.

\vskip 0.1in
\item"{[M]}" R. Mazzeo. The Hodge cohomology of a conformally compact
metric, J. Diff. Geom. 28(1988) 309-339.

\vskip 0.1in
\item"{[Q]}" J. Qing, On the rigidity for conformally compact Einstein
manifolds,
\newline Preprint 2002.

\vskip 0.1in
\item"{[Wa]}" X. Wang, On conformally compact Einstein manifolds,
Math. Res. Lett. 8(2001), no. 5-6, 671--688.

\vskip 0.1in
\item"{[W]}" E. Witten, Anti-de Sitter space and holography, Adv.
Theo. Math. Phy. 2 (1998) 253-290, hep-th/9802150

\vskip 0.1in
\item"{[WY]}" E. Witten and S.T. Yau, Connectedness of the boundary in
ADS/CFT correspondence, Preprint het-th/9910245.

\endroster

\enddocument